%% file: main.tex
\documentclass[letterpaper, 10 pt, conference]{ieeeconf}  

\IEEEoverridecommandlockouts                              
\usepackage{graphics} 
\usepackage{epsfig} 
\usepackage{times} 
\usepackage{amsmath} 
\usepackage{amssymb}  
\newtheorem{lemma}{Lemma}
\newtheorem{theorem}{Theorem}
\newtheorem{remark}{Remark}
\usepackage{cite}
\usepackage{float}
\usepackage{etoolbox}
\AtBeginEnvironment{align}{\small}
\AtBeginEnvironment{equation}{\small}
\usepackage{subcaption}

\title{\LARGE \bf
Non-Gaussian Uncertainty Estimation via Semi-nonparametric Density}

\author{Aaron R. Liao, Kenshiro Oguri, and Michele D. Carpenter
\thanks{This work was funded by The Charles Stark Draper Laboratory under the Draper Scholar program}
\thanks{A. R. Liao is a Draper Scholar and Ph.D. Student, School of Aeronautics and Astronautics, Purdue University, West Lafayette, Indiana, 47907, USA}
\thanks{K. Oguri is Assistant Professor, at the School of Aeronautics and Astronautics, Purdue University, West Lafayette, Indiana, 47907, USA}%
\thanks{M. D. Carpenter is Group Leader and Distinguished Member of the Technical Staff, GNC System Architecture, The Charles Stark Draper
Laboratory, Inc., 555 Technology Square Cambridge, MA 02139.}
}

\begin{document}

\maketitle
\thispagestyle{empty}
\pagestyle{empty}

\begin{abstract}
Accurate and efficient representation of non-Gaussian state distributions is critical for estimation, control, and decision-making, especially when the system is highly nonlinear. This paper develops a mathematical framework for estimating non-Gaussian uncertainty via semi-nonparametric (SNP), or Gallant–Nychka, densities. SNP densities use squared probabilists’ Hermite polynomials to guarantee nonnegative probability densities while representing non-Gaussian behavior. For uncertainty propagation, once fitted to propagated state samples, SNP densities enable direct analytical density and probability evaluations without needing to tune additional hyperparameters besides the polynomial order. To enable their application to nonlinear dynamical systems, we derive analytic marginal and cumulative distribution functions and develop a lossless convex formulation of the coefficient estimation problem along a given sign branch. The approach is demonstrated on the chaotic Lorenz system and compared against Gaussian mixture models and Edgeworth approximations. Results demonstrate that SNP densities can accurately represent non-Gaussian state distributions while enabling direct probability and quantile evaluation, making them a promising tool for density estimation and uncertainty quantification in nonlinear estimation and control.
\end{abstract}

\input{intro}
\input{preliminaries}

\input{methodology_approach}
\input{numerical_results}

\input{conclusion}



\bibliographystyle{IEEEtran}
\bibliography{references_cleaned_IEEE_abbrev}

\end{document}

%% file: intro.tex
\section{Introduction}
Filtering and uncertainty propagation are key elements of estimation, control, and decision-making. A central challenge in these problems is accurately representing and estimating the probability distributions associated with uncertain states \cite{doucet_introduction_2001,kochenderfer_decision_2015}. While density estimation from static datasets is well studied in statistics \cite{silverman_density_1986,scott_multivariate_2015,rosenblatt_remarks_1956,parzen_estimation_1962}, nonlinear dynamical systems introduce additional challenges as initially simple uncertainty distributions can become strongly non-Gaussian. In these settings, uncertainty-quantification methods must accurately represent the evolving state distribution while supporting efficient density and probability evaluations for risk-aware control, probabilistic decision-making, and estimation \cite{liao_higher-order_2026,qi_non-gaussian_2025,kumagai_chance-constrained_2024}.

\par Some existing density estimation methods include classical parametric and nonparametric statistical techniques, Gaussian mixture models (GMMs), particle methods, polynomial chaos expansions (PCEs), and moment-based approximations such as the Edgeworth and Gram--Charlier expansions \cite{silverman_density_1986,scott_multivariate_2015,edgeworth_representation_1907,charlier_uber_1905,vittaldev_multidirectional_2016,doucet_introduction_2001,jones_nonlinear_2013}. Edgeworth and Gram--Charlier expansions provide polynomial corrections to Gaussian densities, but finite-order approximations are not guaranteed to remain nonnegative and therefore may fail to define valid probability densities \cite{barton_conditions_1952,jondeau_gramcharlier_2001}. GMMs and particle methods can represent arbitrary non-Gaussian distributions, but accurately capturing highly nonlinear uncertainty may require many Gaussian kernels or particles, increasing computational complexity \cite{vittaldev_multidirectional_2016,alspach_nonlinear_1972,doucet_introduction_2001}. Adaptive splitting and merging of Gaussian mixands for uncertainty propagation can also add an additional layer of complexity, due to the introduction of additional hyperparameters that must be tuned. PCE provides an attractive method for uncertainty propagation but does not directly yield an explicit probability density without additional effort \cite{jones_nonlinear_2013,xiu2010numerical}. Density evolution may also be computed through the Fokker--Planck equation \cite{risken_fokker-planck_1996}, but this generally requires solving a high-dimensional partial differential equation. These limitations motivate a tractable method for uncertainty quantification that generate valid density estimates while capturing non-Gaussian structure and enabling direct probability evaluation.

    \par Semi-nonparametric (SNP), or Gallant–Nychka, densities are polynomial-based density representations closely related to moment-based approximations such as the Edgeworth expansion. Like Edgeworth expansions, they model departures from Gaussianity using Hermite polynomials; however, the polynomial correction is squared, guaranteeing a nonnegative probability density while retaining the ability to represent non-Gaussian behavior \cite{gallant_semi-nonparametric_1987}. Although SNP densities have been studied extensively in econometrics, their application to uncertainty quantification in nonlinear dynamical systems remains largely unexplored. Such applications require repeated density fitting from propagated state samples together with efficient marginalization and probability evaluation. These capabilities for SNP have not previously been developed in a unified mathematical framework for estimation and control. Once fitted to propagated samples, SNP densities can provide these quantities analytically without the kernel management and tuning procedures associated with mixture-based uncertainty propagation. These properties make SNP densities an attractive alternative option for representing non-Gaussian uncertainty in nonlinear dynamical systems.

\par This paper develops a mathematical, SNP-based framework for density estimation in nonlinear dynamical systems. Rather than proposing a new density family, the contribution of this work is the mathematical and computational framework required to apply SNP densities to nonlinear density estimation and uncertainty propagation. Specifically, the paper develops analytic expressions for marginal SNP probability density functions, an analytic multivariate cumulative distribution function enabling direct probability and quantile evaluation, and a lossless convex formulation of the SNP coefficient estimation problem along a given sign branch. The proposed mathematical framework for uncertainty estimation is demonstrated in the chaotic Lorenz system and quantitatively compared with Gaussian mixture models and Edgeworth density approximations.

%% file: preliminaries.tex
\section{Semi-nonparametric Density Preliminaries}
\subsection{SNP Structure}
The standardized SNP density is given by
\begin{equation}
    p(z)=\frac{\phi(z)P(z)^2}{S},
\end{equation}
where $\phi(z)$ is the standard Gaussian density with zero mean and identity variance, $P(z)$ is a polynomial expansion, and $S=\mathbb{E}_\phi\left[P(z)^2\right]$, is the normalization constant to ensure a valid density function. The subscript $\phi$ on the expectation operator is used to indicate that this expectation is computed with respect to a standard Gaussian random variable. Importantly, since this density is constructed as a perturbation from a standard normal, the input data must also be whitened.
\par The polynomial $P(z)$ is constructed from probabilists' Hermite polynomials, which form an orthogonal basis with respect to the Gaussian weight function. 
The $n$-th order probabilists' Hermite polynomial is defined as
\begin{equation}
H_{n}(z)=(-1)^n e^{z^2/2}\frac{\mathrm{d}^n}{\mathrm{d}z^n}e^{-z^2/2},
\end{equation}
and satisfies the orthogonality condition $\mathbb{E}_\phi\left[H_{m}(Z)H_{n}(Z)\right]=n!\delta_{mn}$. Due to this orthogonality property, Hermite polynomials provide a convenient basis for representing deviations from a Gaussian density, allowing the SNP density to model non-Gaussian distributions while retaining a Gaussian reference.
Unlike the Edgeworth and Gram-Charlier expansions, the added complexity of the SNP density comes in determining the polynomial coefficients, which are no longer functions of the distribution's moments and cumulants.

\subsection{Multivariate SNP Density}
\label{sec:mvarSNP}
The multivariate SNP density is constructed using a set of multi-indices. Let $\alpha\in\mathbb{N}_0^d$ denote a $d$-dimensional multi-index, where $\mathbb{N}_0$ is the set of nonnegative integers. A multi-index is written as $\alpha =
\begin{bmatrix}
\alpha_1 &
\alpha_2 &
\cdots &
\alpha_d
\end{bmatrix}$, and has total degree 
$|\alpha|=\sum_{j=1}^{d}\alpha_j$.
The set of all multi-indices with total degree between $2$ and the truncation order $K$ is

\begin{equation}
\mathcal{A}
=
\left\{
\alpha\in\mathbb{N}_0^d:
2\le |\alpha|\le K
\right\}.
\end{equation}

Since the density is constructed in a whitened coordinate system, the lowest polynomial order is two, as the mean and variance of the distribution are fixed to the zero vector and identity respectively. Enumerating the elements of $\mathcal{A}$ as $\mathcal{A}
=
\left\{
\alpha^{(1)},
\alpha^{(2)},
\dots,
\alpha^{(M)}
\right\}$, the total number of coefficients for a multivariate SNP density of order $K$ and dimension $d$ is $M
=
|\mathcal{A}|
=
\binom{d+K}{K}-1-d$. For convenience, define the reduced multi-index $\alpha_{-k}
=
\begin{bmatrix}
\alpha_1 &
\cdots &
\alpha_{k-1} &
\alpha_{k+1} &
\cdots &
\alpha_d
\end{bmatrix}$, obtained by removing the $k$th component of $\alpha$. We also define the index set $\mathcal{A}_k^{(0)}
=
\left\{
\alpha\in\mathcal{A}:
\alpha_{-k}=\mathbf{0}
\right\}$, and the set of multi-index pairs $\mathcal{P}_k
=
\left\{
(\alpha,\beta)\in\mathcal{A}\times\mathcal{A}:
\alpha_{-k}=\beta_{-k}
\right\}$. For a multi-index $\alpha\in\mathcal{A}$, the corresponding multivariate Hermite polynomial is

\begin{equation}
H_\alpha(\mathbf{z})
=
\prod_{j=1}^{d}
H_{\alpha_j}\!\left((\mathbf{z})_j\right),
\label{eq:mvarHermite}
\end{equation}

where $H_k(\cdot)$ denotes the probabilists' Hermite polynomial of order $k$. Enumerating the elements of $\mathcal{A}$, the coefficient vector is

\begin{equation}
\Theta
=
\begin{bmatrix}
c_{\alpha^{(1)}} &
c_{\alpha^{(2)}} &
\cdots &
c_{\alpha^{(M)}}
\end{bmatrix}^{\top}
\in\mathbb{R}^{M},
\label{eq:mvarTheta}
\end{equation}

and the multivariate Hermite basis vector is

\begin{equation}
\mathcal{H}(\mathbf{z})
=
\begin{bmatrix}
H_{\alpha^{(1)}}(\mathbf{z}) &
H_{\alpha^{(2)}}(\mathbf{z}) &
\cdots &
H_{\alpha^{(M)}}(\mathbf{z})
\end{bmatrix}^{\top}
\in\mathbb{R}^{M}.
\label{eq:mvarH}
\end{equation}

The factorial of a multi-index is defined as $\alpha!
=
\prod_{j=1}^{d}\alpha_j!$, and the normalization matrix is

\begin{equation}
\mathcal{Q}
=
\operatorname{diag}
\left(
\alpha^{(1)}!,
\alpha^{(2)}!,
\ldots,
\alpha^{(M)}!
\right)
\in\mathbb{R}^{M\times M},
\label{eq:mvarQ}
\end{equation}

where each diagonal element is the factorial of the corresponding multi-index in $\mathcal{A}$.
The multivariate polynomial and normalization factor are then

\begin{equation}
P(\mathbf{z};\Theta)
=
1+\Theta^\top\mathcal{H}(\mathbf{z}),
\label{eq:mvarP}
\end{equation}

and

\begin{equation}
S(\Theta)
=
1+\sum_{\alpha\in\mathcal{A}}
c_\alpha^2
\left(
\prod_{i=1}^{d}\alpha_i!
\right)
=
1+\Theta^\top\mathcal{Q}\Theta.
\label{eq:mvarS}
\end{equation}

%% file: methodology_approach.tex
\section{SNP Estimation Framework}

While Gallant--Nychka SNP densities have been studied extensively in econometrics \cite{gallant_semi-nonparametric_1987,gallant_seminonparametric_1989}, relatively little work has considered the fully coupled multivariate formulation required for uncertainty quantification in nonlinear dynamical systems. Existing multivariate extensions, such as the positive Edgeworth--Sargan formulation\cite{niguez_multivariate_2016}, eliminate mixed Hermite products and therefore cannot capture higher-order dependence between variables. The remainder of this section develops the mathematical framework required to apply fully coupled Gallant--Nychka SNP densities to nonlinear uncertainty propagation, including analytic marginal and cumulative distribution functions together with practical coefficient estimation.

\subsection{Marginal SNP Distributions}
\label{sec:marginalSNP}

The multivariate SNP density can be marginalized analytically using the orthogonality of the probabilists' Hermite polynomials. The whitened marginal density in dimension $k$ is
\begin{equation}
p_k(\mathbf{z}_k;\Theta)
=
\frac{1}{S(\Theta)}
\int
\phi_d(\mathbf{z})P(\mathbf{z};\Theta)^2
\prod_{j\neq k}\mathrm{d}\mathbf{z}_j.
\label{eq:marginalPDFog}
\end{equation}

Expanding the squared polynomial,
\begin{equation}
\begin{aligned}
P(\mathbf{z};\Theta)^2
&=
1
+
2\sum_{\alpha\in\mathcal{A}}
c_\alpha H_\alpha(\mathbf{z})
\\
&\quad+
\sum_{\alpha\in\mathcal{A}}
\sum_{\beta\in\mathcal{A}}
c_\alpha c_\beta
H_\alpha(\mathbf{z})H_\beta(\mathbf{z}),
\end{aligned}
\label{eq:squared_polynomials}
\end{equation}
and using independence among the whitened coordinates and Hermite orthogonality, the linear terms survive only when
$\alpha_{-k}=\mathbf{0}$, while the quadratic terms survive only when
$\alpha_{-k}=\beta_{-k}$. Therefore,
\begin{equation}
\begin{aligned}
&p_k(\mathbf{z}_k;\Theta)
=
\frac{\phi(\mathbf{z}_k)}{S(\Theta)}
\Bigg[
1
+
2\sum_{\alpha\in\mathcal{A}_k^{(0)}}
c_\alpha H_{\alpha_k}(\mathbf{z}_k)
\\
&\qquad+
\sum_{(\alpha,\beta)\in\mathcal{P}_k}
c_\alpha c_\beta
\left(
\prod_{j\neq k}\alpha_j!
\right)
H_{\alpha_k}(\mathbf{z}_k)
H_{\beta_k}(\mathbf{z}_k)
\Bigg],
\end{aligned}
\label{eq:marginalSNP}
\end{equation}
where $\mathcal{A}_k^{(0)}$ and $\mathcal{P}_k$ are defined in
Section~\ref{sec:mvarSNP}.

\subsection{Cumulative Distribution Function}
\label{sec:CDF}

Due to the properties of the probabilists' Hermite polynomials, the cumulative distribution function (CDF) in the whitened space of the SNP density can also be derived analytically. The multivariate CDF is defined as $F_Z(\mathbf{z};\Theta)
=
\int_{-\infty}^{z_1}
\cdots
\int_{-\infty}^{z_d}
p(t;\Theta)
\prod_{i=1}^{d}\mathrm{d}t_i$. 
Substituting the SNP density and expanding $P(t;\Theta)^2$ using
\eqref{eq:squared_polynomials} gives
\begin{equation}
F_Z(\mathbf{z};\Theta)
=
\frac{1}{S(\Theta)}
\left(
I_0(\mathbf{z})+I_1(\mathbf{z})+I_2(\mathbf{z})
\right),
\label{eq:CDFdecomposition}
\end{equation}
where
\begin{align}
I_0(\mathbf{z})
&=
\int_{-\infty}^{z_1}
..
\int_{-\infty}^{z_d}
\phi_d(t)
\prod_{i=1}^{d}\mathrm{d}t_i,
\\
I_1(\mathbf{z})
&=
2
\sum_{\alpha\in\mathcal A}
c_\alpha
\int_{-\infty}^{z_1}
..
\int_{-\infty}^{z_d}
\phi_d(t)
H_\alpha(t)
\prod_{i=1}^{d}\mathrm{d}t_i,
\\
I_2(\mathbf{z})
&=
\sum_{\alpha,\beta\in\mathcal A}
c_\alpha c_\beta
\int_{-\infty}^{z_1}
..
\int_{-\infty}^{z_d}
\phi_d(t)
H_\alpha(t)H_\beta(t)
\prod_{i=1}^{d}\mathrm{d}t_i.
\end{align}

Since the standard Gaussian density factorizes as $\phi_d(t)=\prod_{i=1}^{d}\phi(t_i)$, each multidimensional integral separates into a product of one-dimensional integrals. The first integral becomes $I_0(\mathbf{z})
=
\prod_{i=1}^{d}\Phi(z_i)$. Similarly, $I_1(\mathbf{z})
=
2
\sum_{\alpha\in\mathcal A}
c_\alpha
\prod_{i=1}^{d}
\left(
\int_{-\infty}^{z_i}
\phi(t_i)
H_{\alpha_i}(t_i)
\mathrm{d}t_i
\right)$. Using the identity
\begin{equation}
\int_{-\infty}^{x}
\phi(t)H_n(t)\mathrm{d}t
=
\begin{cases}
\Phi(x), & n=0,\\
-H_{n-1}(x)\phi(x), & n\ge1,
\end{cases}
\label{eq:integralPropertyHermite}
\end{equation}
gives $I_1(z)
=
2
\sum_{\alpha\in\mathcal A}
c_\alpha
\prod_{i=1}^{d}
G_{\alpha_i}(z_i)$, where
\begin{equation}
G_n(x)
=
\begin{cases}
\Phi(x), & n=0,\\
-H_{n-1}(x)\phi(x), & n\ge1.
\end{cases}
\end{equation}

For the quadratic term, we use the Hermite product identity
\begin{equation}
H_p(x)H_q(x)
=
\sum_{r=0}^{\min(p,q)}
r!
\binom{p}{r}
\binom{q}{r}
H_{p+q-2r}(x).
\label{eq:HermiteProduct}
\end{equation}

Applying this identity independently in each coordinate yields $I_2(\mathbf{z})
=
\sum_{\alpha,\beta\in\mathcal A}
c_\alpha c_\beta
\prod_{i=1}^{d}
J_{\alpha_i,\beta_i}(z_i)$, where
\begin{equation}
J_{p,q}(x)
=
\sum_{r=0}^{\min(p,q)}
r!
\binom{p}{r}
\binom{q}{r}
G_{p+q-2r}(x).
\end{equation}

Substituting the expressions for $I_0$, $I_1$, and $I_2$ into
\eqref{eq:CDFdecomposition} gives
\begin{equation}
\begin{aligned}
F_Z(\mathbf{z};\Theta)
&=
\frac{1}{S(\Theta)}
\Bigg(
\prod_{i=1}^{d}\Phi(z_i)
+
2
\sum_{\alpha\in\mathcal A}
c_\alpha
\prod_{i=1}^{d}
G_{\alpha_i}(z_i)
\\
&\qquad+
\sum_{\alpha,\beta\in\mathcal A}
c_\alpha c_\beta
\prod_{i=1}^{d}
J_{\alpha_i,\beta_i}(z_i)
\Bigg).
\end{aligned}
\label{eq:mvarSNPCDF}
\end{equation}

\subsection{SNP Coefficient Estimation}
\label{sec:coefficientEstimation}

The SNP-based estimation framework utilizes sample points as data to estimate the underlying density. These points are first whitened according to
\begin{equation}
\mathbf{z}_i=L^{-1}\left(\mathbf{x}_i-\mu\right),
\label{eq:whitePoint}
\end{equation}
where $\mathbf{x}_i$ is a sample point, $\mathbf{z}_i$ denotes the $i$-th whitened sample point, $\mu$ is the sample mean computed from the sample points, and $L$ is the Cholesky square root factor of the sample covariance. Probabilistic evaluations are subsequently performed in the base, whitened space.

Once the sample points are whitened, the remaining task is to estimate the SNP polynomial coefficients. Unlike the Edgeworth and Gram-Charlier expansions, the coefficients in the SNP density are not available in closed-form functions of the cumulants or moments of the random variable. Instead, the SNP coefficients are obtained by maximizing the expected log-likelihood of the SNP density,
\begin{equation}
\hat{\Theta}
=
\arg\max_\Theta
\mathbb{E}\left[
\log\left(
\frac{
\phi(\mathbf{z})P(\mathbf{z};\Theta)^2
}{
\mathbb{E}_\phi\left[P(\mathbf{z};\Theta)^2\right]
}
\right)
\right],
\label{eq:snpMLE}
\end{equation}
where $\Theta$ is defined in \eqref{eq:mvarTheta}. These expectation integrals can be approximated using MC sampling methods \cite{robert_monte_2004,caffo_ascent-based_2005}. MC sampling is employed in this work due to its simplicity and broad applicability. However, alternative sampling strategies such as importance sampling, PCE, cubature, or quadrature may also be employed \cite{owen2013monte,xiu2010numerical}.

Using the weighted and whitened samples, the multivariate SNP coefficient estimation problem can be rewritten as
\begin{equation}
\min_{\Theta\in\mathbb{R}^M}
F(\Theta)
:=
\log\left(1+\Theta^\top Q\Theta\right)
-
2\sum_{i=1}^{N_s}
w_i
\log\left|
1+\mathcal{H}_i^\top\Theta
\right|,
\label{eq:mvarSNPopt}
\end{equation}
where $\mathcal{H}_i=\mathcal{H}(\mathbf{z}_i)$ is the $K$-th order multivariate Hermite polynomial vector defined in \eqref{eq:mvarH} evaluated at the $i$-th whitened sample point and $i=1,\dots,N_s$ where $N_s$ is the total number of samples. The resulting coefficients define the SNP density estimate in the whitened coordinates.

\subsection{Lossless Convexification Along a Sign Branch}
\label{sec:losslessConvexification}

The absolute-value terms in \eqref{eq:mvarSNPopt} separate
the coefficient space into sign branches. Define the sign vector
$\sigma=(\sigma_1,\ldots,\sigma_{N_s})$, $\sigma_i\in\{-1,+1\}$, the signed
affine factors $g_i(\Theta)
    :=
    \sigma_i(1+\mathcal{H}_i^\top\Theta)$, and the corresponding sign branch
\begin{equation}
    \mathcal{D}_\sigma
    :=
    \left\{
    \Theta\in\mathbb{R}^n:
    g_i(\Theta)>0,\; i=1,\ldots,N_s
    \right\}.
    \label{eq:signBranch}
\end{equation}
Each nonempty $\mathcal{D}_\sigma$ is convex since it is the intersection of
affine halfspaces. On a fixed sign branch, the objective becomes
\begin{equation}
    F_\sigma(\Theta)
    =
    \log(1+\Theta^\top Q\Theta)
    -2\sum_{i=1}^{N_s}w_i\log g_i(\Theta).
    \label{eq:branchObjective}
\end{equation}
Using $\sum_i w_i=1$ and the monotonicity of the logarithm, minimizing
\eqref{eq:branchObjective} is equivalent to minimizing
\begin{equation}
    R_\sigma(\Theta)
    :=
    \frac{1+\Theta^\top Q\Theta}
    {\left(\prod_{i=1}^{m}g_i(\Theta)^{w_i}\right)^2},
    \qquad
    \Theta\in\mathcal{D}_\sigma.
    \label{eq:R_sigma}
\end{equation}

\begin{lemma}[Equivalent homogeneous formulation]
\label{lem:homogeneous}
For any fixed, nonempty sign branch $\mathcal{D}_\sigma$, introduce
$t>0$ and $x=t\Theta$, and define
\begin{equation}
    \Gamma_\sigma(x,t)
    :=
    \prod_{i=1}^{m}
    \left[
    \sigma_i(t+\mathcal{H}_i^\top x)
    \right]^{w_i}.
    \label{eq:Gamma}
\end{equation}
Then minimizing \eqref{eq:R_sigma} over $\mathcal{D}_\sigma$ is equivalent to
\begin{subequations}
\label{eq:equalityNormalized}
\begin{align}
    \min_{x,t}\quad
        & t^2+x^\top Qx
        \label{eq:equalityNormalized_obj}\\
    \mathrm{s.t.}\quad
        & \Gamma_\sigma(x,t)=1,
        \label{eq:equalityNormalized_gamma}\\
        & \sigma_i(t+\mathcal{H}_i^\top x)>0,\quad i=1,\ldots,N_s,
        \label{eq:equalityNormalized_sign}\\
        & t>0.
        \label{eq:equalityNormalized_t}
\end{align}
\end{subequations}
The corresponding coefficient vector is recovered as $\Theta=x/t$.
\end{lemma}

\begin{proof}
For $\Theta=x/t$, \begin{equation}
    g_i(x/t)
    =
    \frac{\sigma_i(t+\mathcal{H}_i^\top x)}{t}.
\end{equation} Since $\sum_iw_i=1$, and $\prod_{i=1}^{N_s}g_i(x/t)^{w_i}
    =
    \frac{\Gamma_\sigma(x,t)}{t}$,
then
\begin{equation}
    R_\sigma(x/t)
    =
    \frac{t^2+x^\top Qx}
    {\Gamma_\sigma(x,t)^2}.
    \label{eq:homogeneousRatio}
\end{equation}
The right-hand side is homogeneous of degree zero under
$(x,t)\mapsto(\alpha x,\alpha t)$, $\alpha>0$, while
$\Theta=x/t$ is unchanged. Moreover, $\Gamma_\sigma$ is homogeneous of
degree one. Thus, for any $(x,t)$ with $\Gamma_\sigma(x,t)>0$, choosing
$\alpha=1/\Gamma_\sigma(x,t)$ uniquely scales $(x,t)$ such that
$\Gamma_\sigma(\alpha x,\alpha t)=1$.

Specifically, for any $\Theta\in\mathcal{D}_\sigma$, choose $\hat{t}=\left(\prod_{i=1}^{N_s}g_i(\Theta)^{w_i}\right)^{-1}$, and $\hat{x}=\hat{t}\Theta$. Then $\Gamma_\sigma(\hat{x},\hat{t})=1$ and $\sigma_i(\hat{t}+\mathcal{H}_i^\top\hat{x})= \hat{t}g_i(\Theta)>0$, so the sign branch is preserved. Conversely, any feasible $(x,t)$ in
\eqref{eq:equalityNormalized} yields $\Theta=x/t\in\mathcal{D}_\sigma$.
Using $\Gamma_\sigma(x,t)=1$ in \eqref{eq:homogeneousRatio} establishes
equality of the objective values in both directions.
\end{proof}

The equality constraint in \eqref{eq:equalityNormalized} is
nonconvex. The following result shows that it can be relaxed to a convex
inequality without changing the optimum.

\begin{lemma}[Lossless relaxation]
\label{lem:lossless}
Problem \eqref{eq:equalityNormalized} is equivalent in optimal value to
\begin{subequations}
\label{eq:losslesscvxRelaxProb}
\begin{align}
    \min_{x,t}\quad
        & t^2+x^\top Qx
        \label{eq:losslesscvxRelaxProb_obj}\\
    \mathrm{s.t.}\quad
        & \Gamma_\sigma(x,t)\geq 1,
        \label{eq:losslesscvxRelaxProb_gamma}\\
        & \sigma_i(t+\mathcal{H}_i^\top x)>0,\quad i=1,\ldots,N_s,
        \label{eq:losslesscvxRelaxProb_sign}\\
        & t>0.
        \label{eq:losslesscvxRelaxProb_t}
\end{align}
\end{subequations}
Furthermore, every minimizer of \eqref{eq:losslesscvxRelaxProb} satisfies
$\Gamma_\sigma(x,t)=1$.
\end{lemma}

\begin{proof}
Suppose an optimal solution $(x^\star,t^\star)$ satisfies
$\Gamma_\sigma(x^\star,t^\star)>1$. Define $\alpha
    =
    \Gamma_\sigma(x^\star,t^\star)^{-1}\in(0,1)$ and $(\tilde{x},\tilde{t})
    =
    (\alpha x^\star,\alpha t^\star)$. By degree-one homogeneity, $\Gamma_\sigma(\tilde{x},\tilde{t})=1$.
Since $\alpha>0$, all positivity constraints in \eqref{eq:losslesscvxRelaxProb_sign} and \eqref{eq:losslesscvxRelaxProb_t} remain satisfied. However, $\tilde{t}^{\,2}+\tilde{x}^\top Q\tilde{x}
    =
    \alpha^2
    \left[
    (t^\star)^2+(x^\star)^\top Qx^\star
    \right]
    <
    (t^\star)^2+(x^\star)^\top Qx^\star$, contradicting optimality. Hence every minimizer satisfies
$\Gamma_\sigma(x^\star,t^\star)=1$, and the relaxation is lossless.
\end{proof}

\begin{theorem}[Branchwise convexity]
\label{thm:branchwiseConvexity}
For every fixed, nonempty sign branch $\mathcal{D}_\sigma$,
\eqref{eq:losslesscvxRelaxProb} is a convex optimization problem and is
equivalent in optimal value to the original SNP coefficient estimation
problem restricted to $\mathcal{D}_\sigma$.
\end{theorem}

\begin{proof}
Since $Q\succ0$, the objective $t^2+x^\top Qx$ is strictly convex.
Furthermore, $\Gamma_\sigma(x,t)$ is a weighted geometric mean of the
positive affine functions
$\sigma_i(t+H_i^\top x)$. For $w_i>0$ and $\sum_iw_i=1$, the weighted
geometric mean is concave on the positive orthant. Consequently, $\left\{(x,t):\Gamma_\sigma(x,t)\geq1\right\}$ is convex. The remaining constraints are affine halfspaces. Convexity
therefore follows, while equivalence to the original branch problem follows
from Lemmas~\ref{lem:homogeneous} and~\ref{lem:lossless}.
\end{proof}

Therefore, once a sign branch $\sigma$ is specified, the corresponding SNP
coefficient estimation problem can be solved globally by a single convex
program, with the original coefficients recovered through
$\Theta^\star=x^\star/t^\star$. The full coefficient estimation problem
is then a union over feasible sign branches. Note that this means that the branchwise relaxed convex problem derived here is therefore not a globally optimal solution as it is sign branch dependent.

\begin{remark}[Challenges with global optimization]
Finding the global optimum by enumeration suffers from a severe curse of dimensionality, requiring a convex problem for each feasible branch among $2^{N_s}$ possible sign patterns. With the heuristic $N_s=10^n$, where $n$ is the state dimension, even $n=3$ yields $2^{1000}\approx10^{301}$ possible patterns. As an empirical observation, searches over $10^4$ randomly sampled sign patterns across several distributions found only infeasible branches, illustrating the difficulty of identifying candidate branches.
\end{remark}

\begin{remark}[Selection of all-positive sign branch]
The all-positive branch $\sigma^+$ is always feasible since $\Theta=0$ satisfies $1+\mathcal{H}_i^\top\Theta>0$ for all $i$. Tests of all single-sign-flip branches from $\sigma^+$ also found $\sigma^+$ to be among the best-performing branches for the unimodal, moderately non-Gaussian distributions considered. Based on these empirical observations, we hypothesize that the all-positive branch $\sigma^+$ provides accurate density estimates and use $\sigma^+$ as the sign branch for all SNP fits. Efficient selection of optimal branches, particularly for multimodal distributions, remains a direction for future work.
\end{remark}

%% file: numerical_results.tex
\section{Numerical Examples}
The numerical experiments presented in this section use the Lorenz system $\dot{x}=s(y-x),\;
    \dot{y}=x(\rho-z)-y,\;
    \dot{z}=xy-\beta z,$ with $s=10$, $\rho=28$, $\beta=8/3$. Monte Carlo (MC) samples from
$\mathcal{N}([1,1,1]^\top,\operatorname{diag}(0.3^2,0.3^2,0.3^2))$
are propagated to $T=1.1$.

\subsection{Numerical Validation of Lossless Convexification}
To numerically verify the lossless convexification in Section
\ref{sec:losslessConvexification}, the NLP solution of \eqref{eq:snpMLE},
initialized heuristically, is compared with the convex solution on its
corresponding sign branch across different integration times $T$ and SNP
orders $K$. For each case, the converged NLP sign pattern is used to solve the
corresponding branchwise convex problem.

\begin{figure}[H]
    \begin{subfigure}[b]{1\columnwidth}
        \centering
        \includegraphics[width = 1\columnwidth]{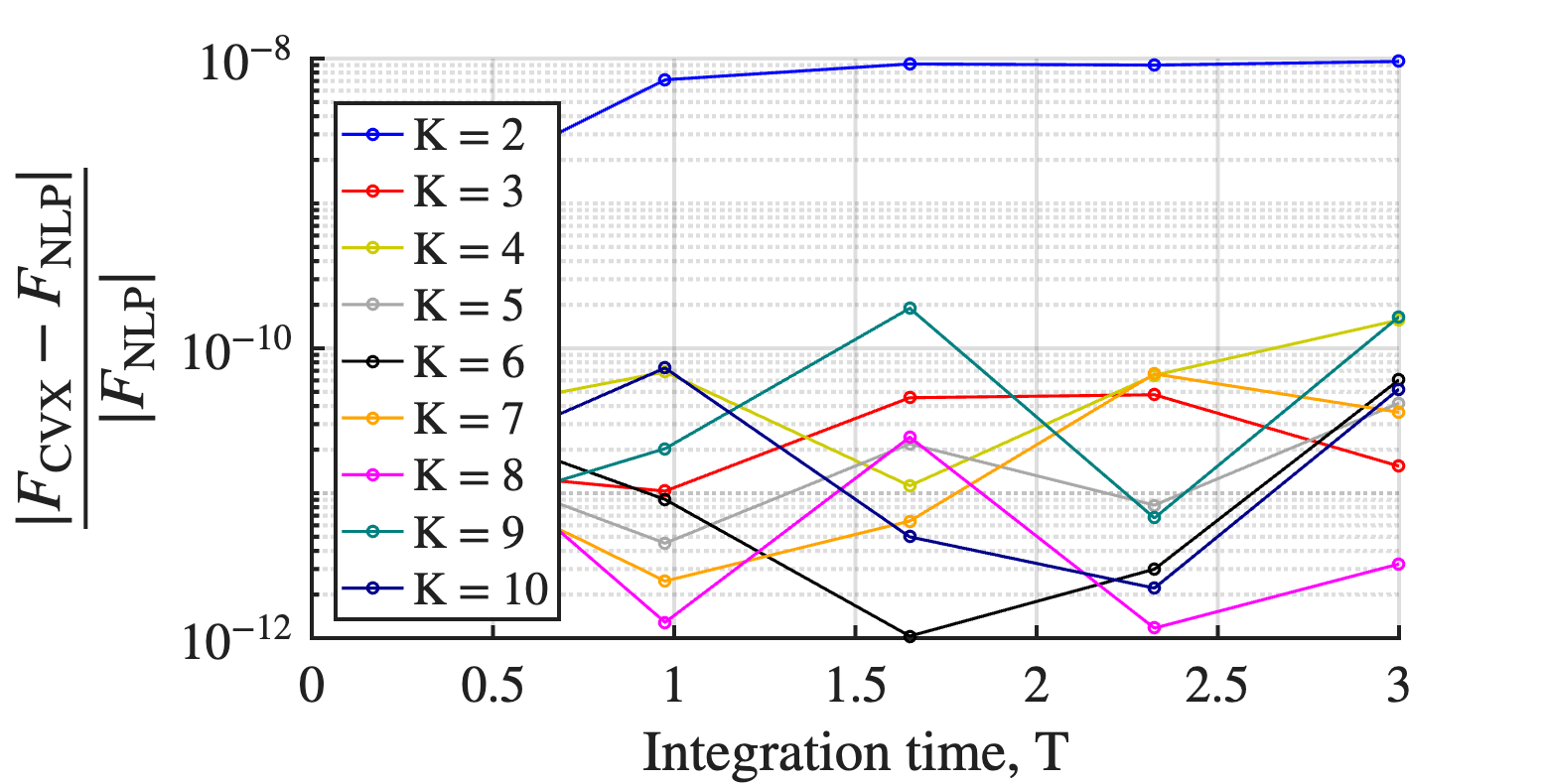}
        \caption{NLP vs.\ convex objective difference}
        \label{fig:NLPvsLosslessDiff}
    \end{subfigure}
    ~
    \begin{subfigure}[b]{1\columnwidth}
        \centering
        \includegraphics[width=1\linewidth]{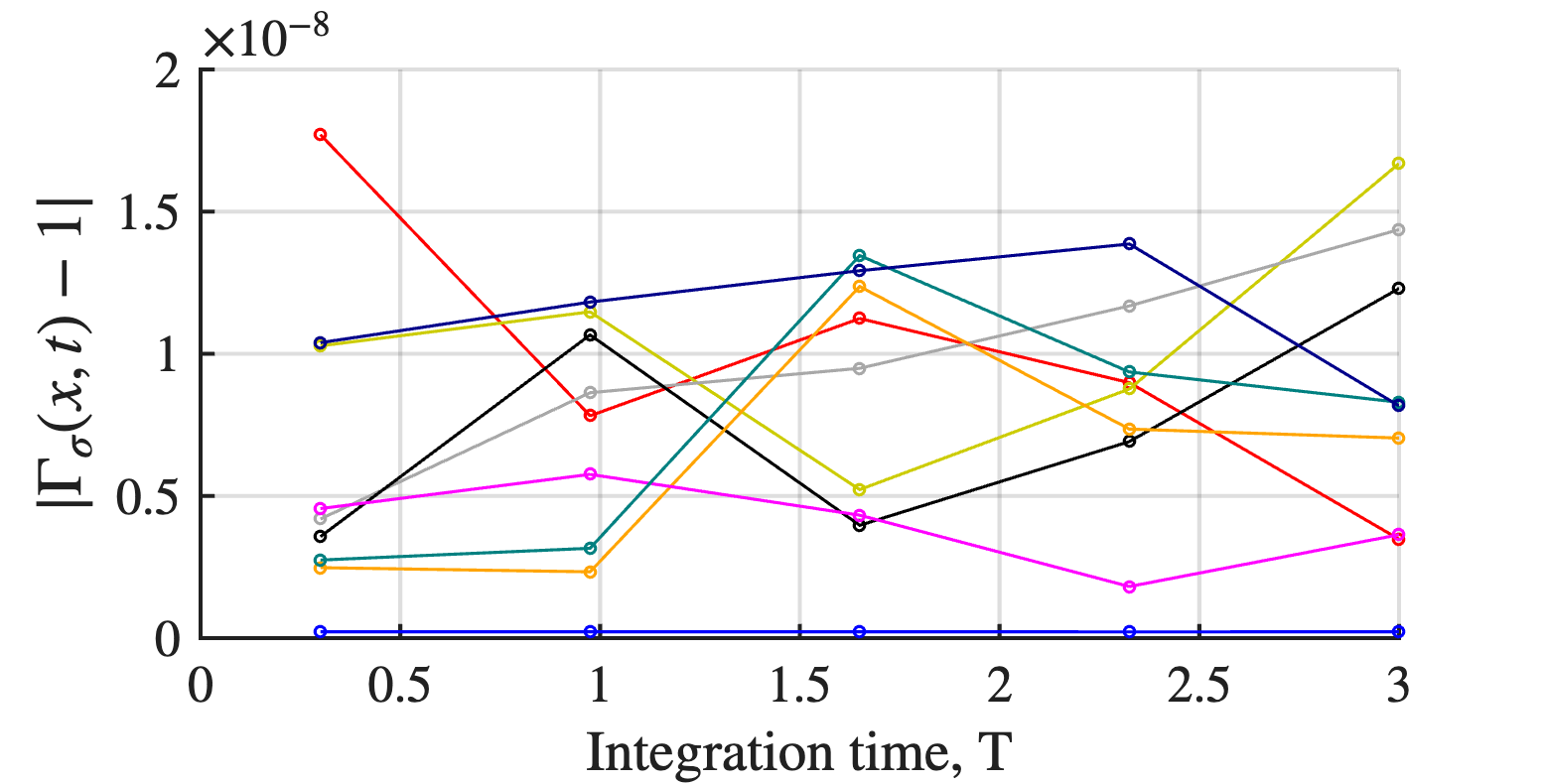}
        \caption{Normalization constraint}
        \label{fig:losslessNormalization}
    \end{subfigure}
    \caption{Numerical validation of losslessness}
    \label{fig:SNPvsEdgeComparisons}
\end{figure}
Fig. \ref{fig:NLPvsLosslessDiff} shows that the normalized objective difference between the NLP and convex solutions is on the order of the solver tolerance, while Fig. \ref{fig:losslessNormalization} verifies that the normalization constraint is active. Together, these results numerically verify the lossless reformulation. Since the reformulated problem is convex, its solution is globally optimal on the specified sign branch.

\subsection{Density Estimation in the Lorenz System}
\label{sec:densityExample}
\begin{figure}[H]
    \centering
    \includegraphics[width=1\columnwidth]{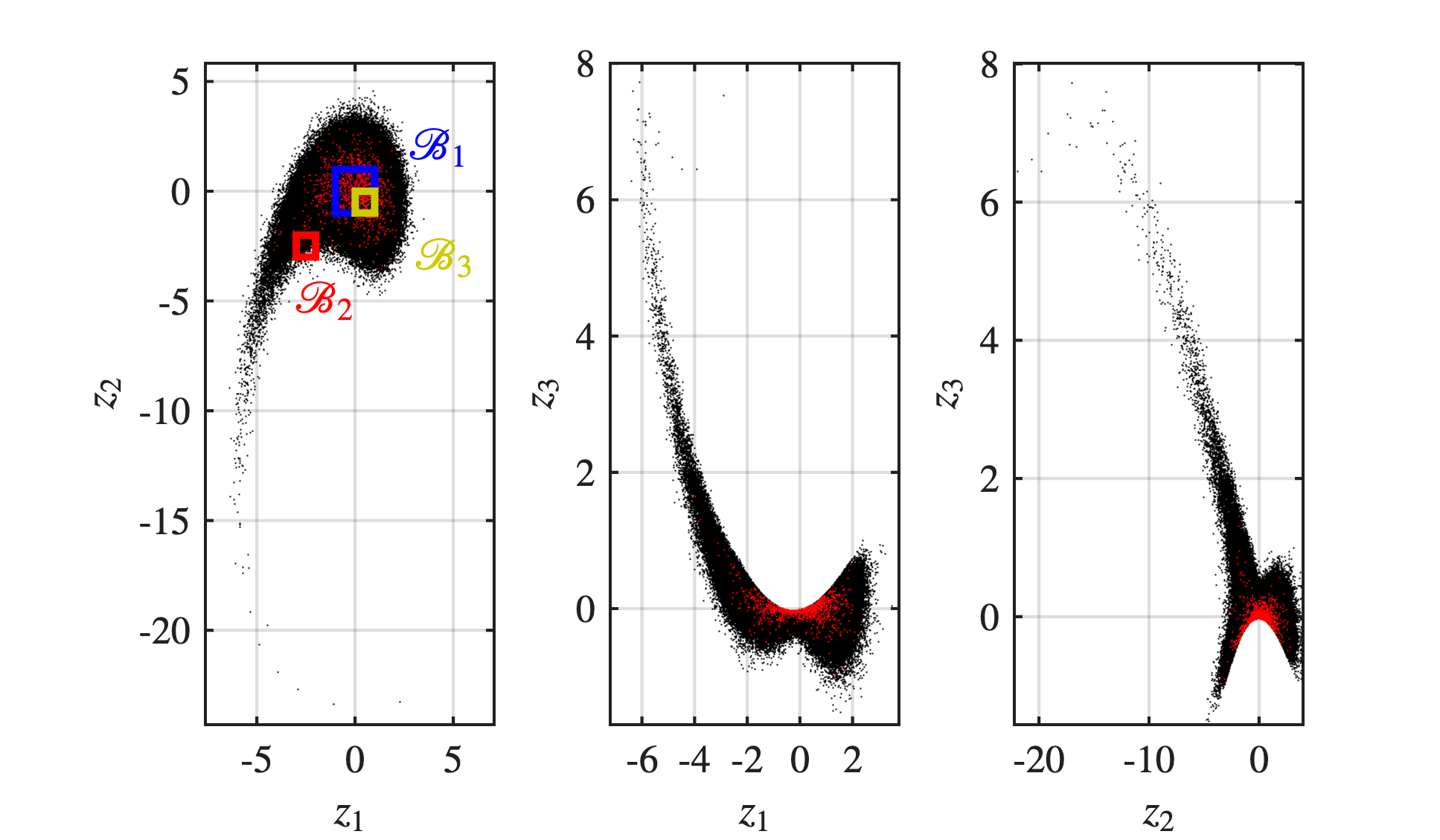}
    \caption{$10^6$ MC Cloud in whitened space}
    \label{fig:mcCloud}
\end{figure}

The $10^6$-sample MC cloud in Fig.~\ref{fig:mcCloud} serves as the
non-Gaussian reference distribution. The red points show an example fitting subset. The boxes, shown in the $z_1$–$z_2$ projection, indicate the regions over which probability mass is evaluated in Section~\ref{sec:CDFeval}.
To assess sensitivity to sample randomness when only a moderate number of propagated samples are available, each density estimator is independently fit to 100 subsets of $N_s=10^3$ samples drawn from the reference MC distribution, consistent with the $10^n$ sampling heuristic for an $n$-dimensional state. The multivariate Edgeworth series are truncated at the 4th order and the coefficients are computed using cumulant tensors computed from the samples. Similarly, GMMs with $N_c$ components and full, component-specific covariances are fit to each sample subset using MATLAB’s \texttt{fitgmdist} implementation of expectation-maximization. For each fit, 10 independent (k)-means++ initializations are performed and the converged solution with the highest-likelihood is retained to reduce sensitivity to local optima. All fitted densities are then evaluated against the full $10^6$ sample reference distribution.

\par The SNP density estimate solved along $\sigma^+$ is quantitatively compared against the multivariate Edgeworth and GMM approximations using the Kullback-Leibler divergence (KLD) between the reference MC cloud to the different density estimates. The KLD is computed as \(D_{\mathrm{KL}}(p_{\mathrm{MC}}\|q)\), where \(p_{\mathrm{MC}}\) denotes the unknown propagated-state density represented by the reference MC cloud and \(q\) is the fitted SNP, Edgeworth, or GMM approximation. After transforming the samples into whitened coordinates, each fitted model is evaluated at every Monte Carlo point to estimate its cross-entropy $\widehat{H}(p_{\mathrm{MC}},q)=-\frac{1}{N_{\mathrm{MC}}}\sum_{i=1}^{N_{\mathrm{MC}}}\log q(z_i)$. The entropy of the reference distribution, \(\widehat{H}(p_{\mathrm{MC}})\), is estimated using a nearest-neighbor entropy estimator and is shared across all models. The divergence is then computed as $\widehat{D}_{\mathrm{KL}}(p_{\mathrm{MC}}\|q)=\widehat{H}(p_{\mathrm{MC}},q)-\widehat{H}(p_{\mathrm{MC}})$. For Edgeworth expansions, nonpositive density values are replaced by a small numerical floor of $10^{-100}$ when computing logarithms.

\par The KLD values for these density estimation methods over the 100 sample sets are tabulated in Table~\ref{tab:kldUQ} below. 
\begin{table}[H]
    \centering
    \caption{Kullback-Leibler divergence statistics}
    \begin{tabular}{c|c|c}
        \textbf{Model} & \textbf{Parameters} & \textbf{KLD (mean $\pm$ std.)} \\
        \hline
        GMM $N_c=2$  & $19$  & $20.132 \pm 8.7225$ \\
        GMM $N_c=4$  & $39$  & $35.423 \pm 22.01$ \\
        GMM $N_c=8$  & $79$  & $93.734 \pm 65.05$ \\
        Edgeworth    & $25$  & $69.726 \pm 69.37$ \\
        SNP $K=4$    & $31$  & $3.2187\pm0.00181$ \\
        SNP $K=6$    & $80$  & $3.057 \pm 0.0027265$ \\
        SNP $K=8$    & $161$ & $2.9319 \pm 0.0037701$ \\
    \end{tabular}
    \label{tab:kldUQ}
\end{table}
Table~\ref{tab:kldUQ} shows that, for the $N_s=10^3$ sample sets, SNP achieves substantially lower mean KLD and trial-to-trial variability than GMM and Edgeworth. In particular, the $K=6$ SNP outperforms the comparable $N_c=8$ GMM. Increasing GMM model complexity increases both the mean KLD and its variability, indicating greater sensitivity to the limited fitting samples. As a consistency check, increasing the fitting-set size to $N_s=10^4$ substantially improves the GMM estimates, with higher-order GMMs beginning to outperform SNP. This supports that the degradation in Table~\ref{tab:kldUQ} is primarily due to sample sensitivity. In contrast to this behavior, increasing SNP order improves the mean density approximation while maintaining low variability.

\subsection{Cumulative Distribution Function Evaluation}
\label{sec:CDFeval}
\begin{figure}[H]
    \centering
    \includegraphics[width=0.96\columnwidth]{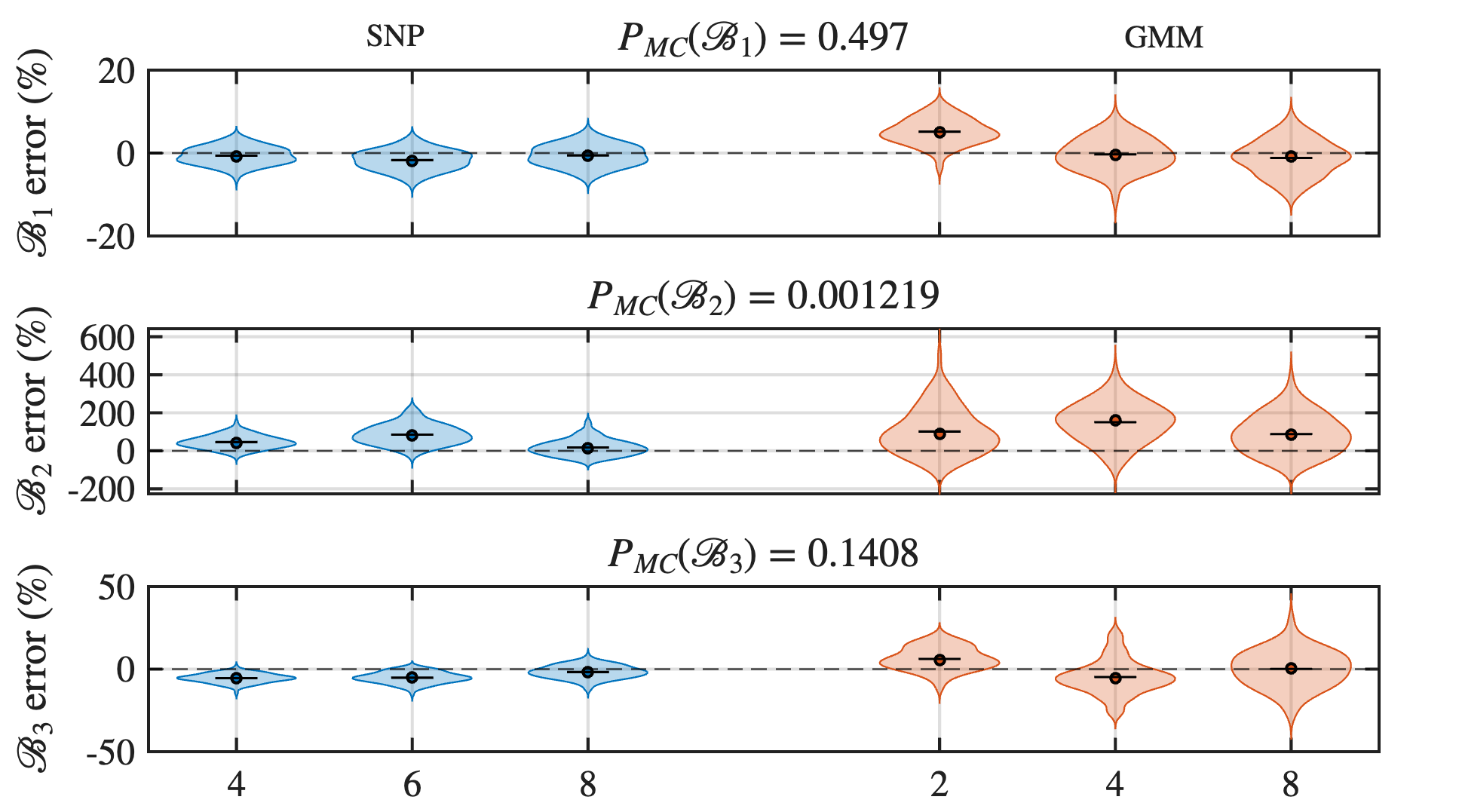}
    \caption{CDF box evaluation violin plots}
    \label{fig:cdfBoxViolinPlots}
\end{figure}

The analytic SNP CDF derived in \eqref{eq:mvarSNPCDF} is evaluated over three regions in the whitened space:
$\mathcal{B}_1:(z_1,z_2)\in[-1,1]\times[-1,1]$,
$\mathcal{B}_2:(z_1,z_2)\in[-3,-2]\times[-3,-2]$, and
$\mathcal{B}_3:(z_1,z_2)\in[0,1]\times[-1,0]$, containing approximately
$49.7\%$, $0.12\%$, and $14\%$ of the reference MC cloud, respectively.
Fig.~\ref{fig:cdfBoxViolinPlots} shows the resulting distribution of the probability estimate error defined as $\left(P_\mathrm{estimate}(\mathcal{B})-P_\mathrm{MC}(\mathcal{B})\right)/P_\mathrm{MC}(\mathcal{B})$
across the trials.
The SNP estimates exhibit less trial-to-trial variability than the GMM
estimates despite similar mean errors. Table~\ref{tab:cdfUQ} reports the
probability RMSE defined as $100*\mathrm{RMSE}$.
\begin{table}[H]
    \centering
    \caption{CDF box probability RMSE (percentage points)}
    \setlength{\tabcolsep}{3pt}
    \begin{tabular}{c|ccc|ccc|c}
        & \multicolumn{3}{c|}{SNP Order $K$}
        & \multicolumn{3}{c|}{GMM $N_c$} 
        & Edge. \\
        Box 
        & $4$ & $6$ & $8$
        & $2$ & $4$ & $8$
        & \\
        \hline
        $\mathcal{B}_1$
        & $1.08$ & $1.44$ & $1.27$
        & $2.99$ & $1.92$ & $2.12$
        & $16.09$ \\
        $\mathcal{B}_2$
        & $0.069$ & $0.121$ & $0.055$
        & $0.189$ & $0.219$ & $0.161$
        & $1.48$ \\
        $\mathcal{B}_3$
        & $0.864$ & $0.870$ & $0.591$
        & $1.30$ & $1.59$ & $1.53$
        & $12.06$ \\
    \end{tabular}
    \label{tab:cdfUQ}
\end{table}
The SNP and GMM estimates both provide reasonably accurate probability estimates with $N_s=1000$, despite the larger global KLD of the GMM fits. However, the probability RMSEs of the GMM estimates are approximately 2–3 times worse than the SNP estimates. Notably, the SNP estimates maintain small absolute probability errors even for the low-probability region $\mathcal{B}_2$. In contrast, the Edgeworth approximation exhibits substantially larger errors. These results demonstrate that the analytical SNP CDF enables accurate probability evaluations directly from the fitted density.

%% file: conclusion.tex
\section{CONCLUSIONS}
In this paper, we develop a mathematical and computational framework for applying semi-nonparametric (SNP) densities to uncertainty quantification. We derive analytically evaluable marginal densities and multivariate cumulative density functions for SNP densities and show that coefficient estimation can be solved globally on a specified sign branch. The proposed approach represents samples propagated through nonlinear dynamics with an analytical density, enabling direct evaluation of density and probability quantities relevant to estimation and control. Numerical results for the Lorenz system demonstrate accurate and consistent SNP density and probability estimates from a moderate number of samples relative to GMM and Edgeworth approximations. These properties make SNP densities a promising option for uncertainty quantification in nonlinear estimation and control problems requiring valid non-Gaussian density representations and efficient probability evaluation.